\documentclass[english,10pt]{kart8}
\usepackage{babel}
\usepackage{epsfig,graphicx,times}   
\usepackage{amssymb,amsmath}
\def\Ps{\mathcal{P}}
\def\di{m}
\def\Ed{\mathbb{E}^\di}
\newcommand{\diam}{{\mbox{{diam}}}}

\newenvironment{enumeratea}{\begin{list}{(\alph{enumi})}%
{\usecounter{enumi}%
\setlength{\topsep}{0mm}%
\setlength{\partopsep}{0mm}%
\setlength{\itemsep}{0mm}%
\setlength{\labelsep}{2mm}%
\setlength{\labelwidth}{5mm}%
\setlength{\leftmargin}{0mm}%
\addtolength{\leftmargin}{\labelwidth}%
\addtolength{\leftmargin}{\labelsep}%
\setlength{\itemindent}{0mm}%
}}{\end{list}}

\date{May 7, 2019}

\begin{document}
\begin{frontmatter}
  \title{Bounds for the minimum diameter of integral point sets}

  \author{Sascha Kurz and Reinhard Laue}
  \ead{\{sascha.kurz, reinhard.laue\}@uni-bayreuth.de}
  \address{University of Bayreuth, Department of Mathematics, D-95440 Bayreuth, Germany}

  \begin{abstract}
    Geometrical objects with integral sides have attracted mathematicians for ages. For example, the problem 
    to prove or to disprove the existence of a perfect box, that is, a rectangular parallelepiped with all edges, 
    face diagonals and space diagonals of integer lengths, remains open. More generally an integral point set $\Ps$ 
    is a set of $n$ points in the $\di$-dimensional Euclidean space $\Ed$ with pairwise integral distances where the 
    largest occurring distance is called its diameter. From the combinatorial point of view there is a natural interest 
    in the determination of the smallest possible diameter $d(\di,n)$ for given parameters $\di$ and $n$. We give some 
    new upper bounds for the minimum diameter $d(\di,n)$ and some exact values.
  \end{abstract}
  \begin{keyword}
    integral distances \sep diameter  
    \MSC 52C10* \sep 11D99 
  \end{keyword}
\end{frontmatter}

\section{Introduction}

Geometrical objects with integral sides have long attracted mathematicians. One of the earliest results is due to the Pythagoreans and characterizes the smallest rectangle with integral sides and diagonals, more precisely, the integral rectangle with the smallest possible diameter where diameter denotes the largest occurring distance of the points. This is a rectangle with edge lengths $3$ and $4$ so that the diagonal has length $5$ by Pythagoras' Theorem. In this context, a famous old open problem is to show the existence of a perfect box, a rectangular parallelepiped with all edges, face diagonals and space diagonals of integer lengths \cite{UPIN,unsolved_1}. Because this problem seems to be too hard for our current state of mathematics, the authors of \cite{0822.51010} considered combinatorial boxes, i.e., convex bodies with six quadrilaterals as faces, and gave $20$ examples of integral combinatorial boxes, one of which is proven to be minimal with regard to the diameter in \cite{combinatorial_box}. In \cite{parameter_hexaeder}, it is shown that there exist infinitely many integral combinatorial boxes.

Generally, an integral point set $\Ps$ is a set of $n$ points in the $\di$-dimensional Euclidean space $\Ed$ with pairwise integral distances, where not all $n$ points are contained in a hyperplane. From the combinatorial point of view, there is a natural interest in the minimum possible diameter $d(\di,n)$ for given parameters $\di$ and $n$.

In the following, we will focus on bounds and exact numbers for $d(\di,n)$. For a more general overview and applications on integral point sets and similar structures, we refer to \cite{integral_distances_in_point_sets}. Clearly, the condition $n\ge \di+1$ is necessary for an $\di$-dimensional point set. Due to general constructions, see i.e. \cite{diameters}, the condition is also sufficient for the existence of an $\di$-dimensional integral point set consisting of $n$ points.

\begin{thm}
  \label{thm_known_results}
  For $n\ge\di+1$ we have
  \begin{enumeratea}
    \item  $d(\di,n)\le
          \begin{cases}
            2^{n-\di+1}-2 & \mbox{for } n-\di\equiv 0\mod 2,\\
            3(2^{n-\di}-1) & \mbox{for } n-\di\equiv 1\mod 2,
            \end{cases}$\hfill{\cite{diameters}}
    \item $d(\di,n)\le (n-\di)^{c\log\log (n-\di)}$ for a sufficiently large 
          constant $c$, \hfill{\cite{minimum_diameter}}         
    \item $\sqrt{\frac{3}{2\di}}n^{1/\di}<d(\di,n)$, \hfill{\cite{kanold}}
    \item $\frac{1}{\sqrt{14}}n^{1/2}<d(3,n)$ for $n\ge 5$, \hfill{\cite{kanold}}
    \item $cn\le d(2,n)$ for a sufficiently small constant $c$, \hfill{\cite{note_on_integral_distances}}
    \item $d(n,n+1)=1$,
    \item $3\le d(\di,n)\le 4$ for $\di+2\le n\le 2\di$ and $d(\di,2\di)=4$, \hfill{\cite{small_diameter,dipl_piepmeyer}}
    \item $d(\di,2\di+1)\le 8$, \hfill{\cite{dipl_piepmeyer}}
    \item $d(\di,2\di+2)\le 13$, \hfill{\cite{dipl_piepmeyer}}
    \item $d(\di,3\di)\le 109$,  \hfill{\cite{hab_kemnitz}}
    \item and $d(\di,n-1)\le d(\di,n)$.
  \end{enumeratea}
\end{thm}

We conjecture that $d(\di-1,n)\ge d(\di,n)$. Each of the known bounds 
are increasing in $n$ for fixed $\di$ and decreasing in $\di$ for fixed $n$. Several functional relations $f$ between $\di$ and $n$ 
exist for which $d(\di,f(\di))$ can be bounded from above by a constant. Examples are the inequalities of Theorem  \ref{thm_known_results}.(g,h,i,j) and of Theorem \ref{thm_core}.(a) below.

Aside from general bounds, some exact values of $d(\di,n)$ have been determined 
(the bold printed value $d(3,9)=16$ was incorrectly stated as $d(3,9)=17$ in the literature, see i.e. \cite{integral_distances_in_point_sets,dipl_piepmeyer}):

$\left(d(2,n)\right)_{n=3,\dots ,89}= 1,4,7,8,17,21,29,40,51,63,74,91,104,121,134,153,164,$\\
$196,212,228,244,272,288,319,332,364,396,437,464,494,524,553,578,608,$\\
$642,667,692,754,816,897,959,1026,1066,1139,1190,1248,1306,1363,1410,$\\
$1460,1514,1564,1614,1675,1727,1770,1817,1887,1906,2060,2140,2169,$\\
$2231,2299,2432,2494,2556,2624,2692,2827,2895,2993,3098,3196,3294,$\\
$3465,3575,3658,3749,3885,3922,4223,4380,4437,4559,4693,4883$           \hfill{\cite{integral_distances_in_point_sets,phd_kurz,paper_alfred}}

$\left(d(3,n)\right)_{n=4,\dots ,23}=1,3,4,8,13,\mathbf{16},17,17,17,56,65,77,86,99,112,133,154,$\\
$195,212,228$    \hfill{\cite{integral_distances_in_point_sets,paper_characteristic,phd_kurz,dipl_piepmeyer}}

$d(3,5)=d(6,8)=d(8,10)=3$ \hfill{\cite{integral_distances_in_point_sets}}

$d(\di,\di+2)=3$ for $8\le\di\le 23$ \hfill{\cite{phd_kurz}}

$d(\di,n)=4$ for $3\le\di\le 12$ and $\di+3\le n\le 2\di$ \hfill{\cite{phd_kurz}}

$d(\di,n)=4$ for $13\le\di\le 23$ and $2\di-9\le n\le 2\di$ \hfill{\cite{phd_kurz}}

\vspace*{1cm}

Our main results are 

\begin{thm}$\,$\\[-8mm]\\
  \label{thm_core}
  \begin{enumeratea}
    \item $d(\di,\di^2+\di)\le 17$,
    \item $d(\di,n-2+\di)\le\d(2,n)\,\,\,\mbox{for $9\le n\le 122$}$,
  \end{enumeratea}
\end{thm}

the exact values

$d(2,n)_{n=90,\dots,122}=5018,5109,5264,5332,5480,5603,5738,5938,5995,6052,$\\
$6324,6432,6630,6738,6939,7061,7245,7384,7568,7752,7935,8119,8321,$\\
$8406,8648,8729,8927,9052,9211,9423,9534,9794,9905$

$d(3,24)=244$, 

and the following two
\textit{constructions}: 

\begin{thm}
  \label{theo_blow_up_n_minus_1}
  If $\Ps$ is a plane integral point set with diameter $\diam(\Ps)$ consisting of $n$ points, where $n-1$ points
  are situated on a line $\overline{AB}$, then $d(\di,n-2+\di)\le\diam(\Ps)$.
\end{thm}

\begin{thm}
  \label{thm_line_circle}
  If $\Ps$ is a planar integral point set consisting of $n$ points, where $n-1$ points are situated
  on a line $\overline{AB}$, the $n$-th point has distance $h$ to the line $\overline{AB}$, and $\Ps'$
  is an $(\di-1)$-dimensional point set consisting of $n'$ points on an $(\di-1)$-dimensional sphere of radius $h$,
  then we have for $\di\ge 2$ that
  $$
    d(\di,n+n'-1)\le \max(\diam(\Ps),\diam(\Ps')).
  $$ 
\end{thm}

Aside from these results, we have:

{\bf Conjecture}\\[-8mm]
\textit{
  \begin{enumeratea}
    \item $d(\di,n)>(n-\di)^{c\log\log (n-\di)}$ for each fixed $\di$ and suitable large $n$ and $c$,
    \item $d(\di,\di+2)=3$ for $\di\ge 8$,
    \item $d(\di-1,n)\ge d(\di,n)$,
    \item the bound of Theorem \ref{theo_blow_up_n_minus_1} is sharp for $\di=2,\,n\ge 9$; $\di=3,\,n\ge 22$, and 
    $\di\ge 4,\,n\ge \di^2+\di+1$, respectively,
    \item $d(\di,n-2+\di)\le d(2,n)$ for $\di\ge 2$.
  \end{enumeratea}
}


\section{Proofs}
\label{sec_proofs} 

The exact values of $d(\di,n)$ were obtained by exhaustive enumeration via the methods described in \cite{paper_characteristic,phd_kurz,paper_alfred}. For future improvements due to faster computers, we
refer the reader to \cite{own_hp}. By a look at the plane integral point sets with diameter at most $10000$,
it turns out that those with minimum diameter and $9\le n\le 122$ points

\begin{figure}[ht]
  \begin{center}
    \setlength{\unitlength}{0.7cm}
    \begin{picture}(9,4.5)
      \put(0,0){\line(1,0){9}}
      \put(0,0){\line(5,4){5}} 
      \put(1,0){\line(1,1){4}}
      \put(2,0){\line(3,4){3}}
      \put(4,0){\line(1,4){1}}
      \put(6,0){\line(-1,4){1}}
      \put(7,0){\line(-1,2){2}}
      \put(9,0){\line(-1,1){4}}
      \put(0,0){\circle*{0.3}} 
      \put(1,0){\circle*{0.3}} 
      \put(2,0){\circle*{0.3}} 
      \put(4,0){\circle*{0.3}} 
      \put(6,0){\circle*{0.3}} 
      \put(7,0){\circle*{0.3}}  
      \put(9,0){\circle*{0.3}} 
      \put(9.3,-0.2){$l$}
      \put(5.2,4){$P$}
      \put(5,4){\circle*{0.3}} 
      \put(5,0){\line(0,1){4}}
      \put(0.1,-0.5){$a_s$}
      \put(1.1,-0.5){$a_2$}
      \put(2.8,-0.5){$a_1$}
      \put(4.4,0.2){$q$}
      \put(5.2,0.2){$q'$}
      \put(5.1,1.4){$h$}
      \put(4.7,-0.5){$a_0$}
      \put(6.2,-0.5){$a_1'$}
      \put(7.6,-0.5){$a_t'$}
      \put(1.6,2){$b_s$}
      \put(0.8,0.4){$b_2$}
      \put(1.7,0.4){$b_1$}
      \put(3.5,0.4){$b_0$}
      \put(6.05,0.4){$b_0'$}
      \put(7.0,0.4){$b_1'$}
      \put(7,2){$b_t'$}
    \end{picture}\\[2mm]
  \end{center}
  \caption{Plane integral point set $\mathcal{P}$ with $n-1$ points on a line.}
  \label{fig_n_minus_1_points_on_line}
\end{figure}
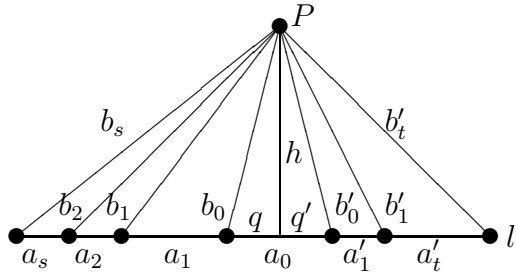

\noindent 
have a very simple structure \cite{phd_kurz,paper_alfred}. They consist of $n-1$ points situated on a line $l$ plus one point $P$ apart from $l$, see Figure \ref{fig_n_minus_1_points_on_line}. An easy method is given in \cite{phd_kurz,paper_alfred} to construct such integral point sets with diameters at most $n^{c\log\log n}$ for a suitably large constant $c$, by choosing integers $h^2$ with many divisors.
%
%
%
If we replace the point $P$ by an $(\di-2)$-dimensional regular simplex
$\mathcal{S}$ with edge length $1$, we obtain an $\di$-dimensional integral point set with the same diameter, which
proves Theorem \ref{theo_blow_up_n_minus_1} and Theorem \ref{thm_core}(b).

\begin{figure}[!ht]
  \begin{center}
    \setlength{\unitlength}{1.0cm}
    \begin{picture}(4,1.5)
      \put(0,0){\circle*{0.2}}
      \put(1,1){\circle*{0.2}}
      \put(3,1){\circle*{0.2}}
      \put(0.46,1.1){$P_1$}
      \put(3.1,1.1){$P_2$}
      \put(4,0){\circle*{0.2}}
      \put(0.8,0){\circle*{0.2}}
      \put(1.5,0){\circle*{0.2}}
      \put(2.7,0){\circle*{0.2}}
      \put(3.2,0){\circle*{0.2}}
      \put(1.86,-0.1){$\dots$}
      \put(2,0.1){$a$}
      \put(0.9,0.39){$b$}
      \put(2.9,0.41){$e$}
      \put(1.9,1.1){$f$}
      \put(0.17,0.4){$d$}
      \put(3.63,0.4){$c$}
      \put(0,0){\line(1,0){4}}
      \put(1,1){\line(1,0){2}}
      \put(0,0){\line(1,1){1}}
      \put(4,0){\line(-1,1){1}}
      \put(0,0){\line(3,1){3}}
      \put(4,0){\line(-3,1){3}}
    \end{picture}
    \caption{Plane point set with points on two parallel lines.}
    \label{fig_two_lines}
  \end{center}
\end{figure}
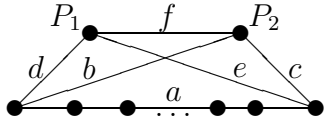

\begin{figure}[ht]
  \begin{center}
  \begin{minipage}{5.5cm}
    \setlength{\unitlength}{0.38cm}
    \begin{picture}(13,11)
      \put(0.000000,4.555764){\circle*{0.4}}
      \put(13.000000,4.555764){\circle*{0.4}}
      \put(11.615385,10.393815){\circle*{0.4}}
      \put(10.943586,10.819940){\circle*{0.4}}
      \put(3.656413,5.838051){\circle*{0.4}}
      \put(5.041028,0.000000){\circle*{0.4}}
      \put(7.958972,9.111528){\circle*{0.4}}
      \put(2.984614,6.264176){\circle*{0.4}}
      %
      \qbezier(0.000000,4.555764)(1.828206,5.196907)(3.656413,5.838051)       
      \qbezier(0.000000,4.555764)(2.520514,2.277882)(5.041028,0.000000)       
      \qbezier(0.000000,4.555764)(1.492307,5.409970)(2.984614,6.264176)       
      \qbezier(13.000000,4.555764)(12.307692,7.474790)(11.615385,10.393815)   
      \qbezier(13.000000,4.555764)(11.971793,7.687852)(10.943586,10.819940)   
      \qbezier[30](13.000000,4.555764)(8.328206,5.196907)(3.656413,5.838051)      
      \qbezier[30](0.000000,4.555764)(6.5,4.555764)(13.000000,4.555764)     
      \qbezier(13.000000,4.555764)(9.020514,2.277882)(5.041028,0.000000)      
      \qbezier(11.615385,10.393815)(11.279485,10.606877)(10.943586,10.819940) 
      \qbezier(11.615385,10.393815)(7.635899,8.115933)(3.656413,5.838051)     
      \qbezier(10.943586,10.819940)(9.451279,9.965734)(7.958972,9.111528)     
      \qbezier(3.656413,5.838051)(4.348720,2.919025)(5.041028,0.000000)       
      \qbezier(7.958972,9.111528)(5.471793,7.687852)(2.984614,6.264176)       
      %
      \put(5.81,3.67){$\mathbf{13}$}  
      \put(2.73,4.8){$\mathbf{9}$}    
      \put(2.62,2.28){$\mathbf{9}$}   
      \put(1.49,5.71){$\mathbf{3}$}   
      \put(12.41,7.47){$\mathbf{6}$}  
      \put(11.17,7.79){$\mathbf{9}$}  
      \put(8.33,5.3){$\mathbf{10}$}   
      \put(9.02,1.58){$\mathbf{8}$}   
      \put(11.28,10.81){$\mathbf{9}$} 
      \put(7.64,7.42){$\mathbf{8}$}   
      \put(9.45,10.37){$\mathbf{3}$}   
      \put(4.45,2.92){$\mathbf{6}$}   
      \put(5.47,8.23){$\mathbf{5}$}   
    \end{picture}
  \end{minipage}
  \begin{minipage}{4.4cm}  
    \setlength{\unitlength}{0.39cm}
      \begin{picture}(11,8.1)
      \put(0,0.4){\circle*{0.55}}
      \put(3,0.4){\circle*{0.55}}
      \put(8,0.4){\circle*{0.55}}
      \put(11,0.4){\circle*{0.55}}
      \put(0,0.4){\line(1,0){11}}
       \put(11,0.4){\qbezier(0,0)(-0.75,3.465)(-1.5,6.93)}
      \put(1.5,7.33){\circle*{0.55}}
      \put(9.5,7.33){\circle*{0.55}}
      \put(1.5,7.33){\line(1,0){8}}
      \put(0,0.4){\qbezier(0,0)(0.75,3.465)(1.5,6.93)}
      \put(3,0.4){\qbezier(0,0)(-0.75,3.465)(-1.5,6.93)}
      \put(8,0.4){\qbezier(0,0)(-3.25,3.465)(-6.5,6.93)}
      \put(11,0.4){\qbezier(0,0)(-4.75,3.465)(-9.5,6.93)}
      \put(0,0.4){\qbezier(0,0)(4.75,3.465)(9.5,6.93)}
      \put(3,0.4){\qbezier(0,0)(3.25,3.465)(6.5,6.93)}
      \put(8,0.4){\qbezier(0,0)(0.75,3.465)(1.5,6.93)}
      \put(1.4,-0.4){$\mathbf{3}$}
      \put(5.4,-0.4){$\mathbf{5}$}
      \put(9.4,-0.4){$\mathbf{3}$}
      \put(5.9,7.4){$\mathbf{8}$}
      \put(0.4,5.1){$\mathbf{9}$}
      \put(2.1,5.1){$\mathbf{9}$}
      \put(2.7,4.21){$\mathbf{11}$}
      \put(3.1,6.3){$\mathbf{13}$}
      \put(0.9,1.9){$\mathbf{13}$}
      \put(3.25,1.9){$\mathbf{11}$}
      \put(8.4,5.1){$\mathbf{9}$}
      \put(10.1,5.1){$\mathbf{9}$}
    \end{picture}  
  \end{minipage}  
  \end{center}
  \caption{$3$-dimensional integral point set consisting of $8$ points with minimum diameter.}
  \label{fig_3_8_min_1}
\end{figure}
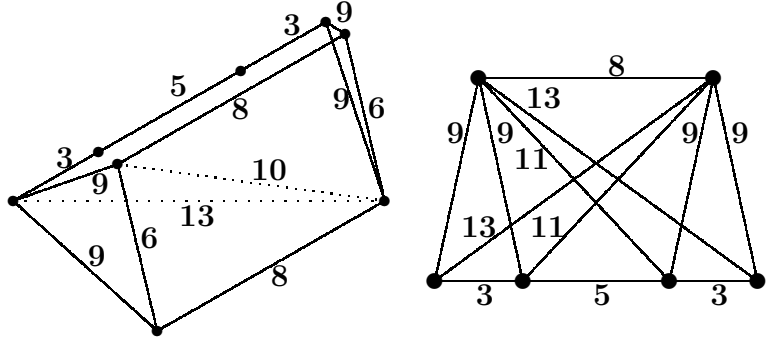

If we assume that we have a plane integral point set $\Ps$ consisting of a line $L$ with $n-2$ points and a parallel line with 
two points $P_1$ and  $P_2$ (see Figure \ref{fig_two_lines}), we can 

\noindent
slightly modify the construction of Theorem \ref{theo_blow_up_n_minus_1} and blow up $P_1$ and $P_2$ to regular $(\di-2)$-dimensional simplices $\mathcal{S}_1$ and $\mathcal{S}_2$ 
of side length $v$. An example is given in Figure \ref{fig_3_8_min_1}. Because the distance of two points $p_1\in\mathcal{S}_1$ and $p_2\in\mathcal{S}_2$ is either $f$ or 
$w:=\sqrt{f^2+v^2}$ we have to choose a suitable $v$ so that $w$ is an integer. 

\begin{thm}
  \label{theo_blow_up_2}
  Let $\Ps$ be a plane integral point set consisting of $n-2$ points on a line $L$ and two points $P_1$ and 
  $P_2$ on a parallel line $M$ with distance $r$ between $L$ and $M$. If there exist positive integers $v,w$ with $f^2+v^2=w^2$ and $v<2r$, 
  where $\overline{P_1P_2}=f$, then $$d(\di,n-2+2(\di-1))\le\max(w,\diam(\Ps)).$$
\end{thm}

Theorem \ref{theo_blow_up_2} is tight in the cases $\di=2$, $n=4,7,8$ and $\di=3$, $n=8$, and also gives $d(\di,2\di+2)\le 13$ (cf. Figure \ref{fig_3_8_min_1}).

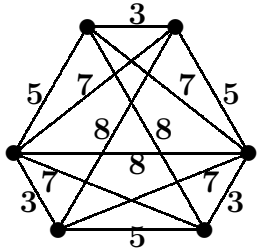
\begin{figure}[ht]
  \begin{center}
    \setlength{\unitlength}{0.39cm}
    \begin{picture}(8,8.1)
      \put(1.5,0.4){\line(1,0){5}}
      \put(6.5,0.4){\circle*{0.55}}
      \put(1.5,0.4){\circle*{0.55}}
      \put(0,3){\circle*{0.55}}
      \put(8,3){\circle*{0.55}}
      \put(1.5,0.4){\qbezier(0,0)(-0.75,1.3)(-1.5,2.6)}
      \put(6.5,0.4){\qbezier(0,0)(0.75,1.3)(1.5,2.6)}
      \put(2.5,7.33){\circle*{0.55}}
      \put(5.5,7.33){\circle*{0.55}}
      \put(2.5,7.33){\line(1,0){3}}
      \put(0,3){\qbezier(0,0)(1.25,2.165)(2.5,4.33)}
      \put(8,3){\qbezier(0,0)(-1.25,2.165)(-2.5,4.33)}
      \put(0,3){\line(1,0){8}}
      \put(3.9,-0.2){$\mathbf{5}$}
      \put(3.9,2.3){$\mathbf{8}$}
      \put(3.9,7.4){$\mathbf{3}$}
      \put(0.2,1){$\mathbf{3}$}
      \put(7.25,1){$\mathbf{3}$}
      \put(0.4,4.7){$\mathbf{5}$}
      \put(7.1,4.7){$\mathbf{5}$}
      \put(0,3){\qbezier(0,0)(2.75,2.165)(5.5,4.33)}
      \put(8,3){\qbezier(0,0)(-2.75,2.165)(-5.5,4.33)}
      \put(1.5,0.4){\qbezier(0,0)(2,3.465)(4,6.93)}
      \put(6.5,0.4){\qbezier(0,0)(-2,3.465)(-4,6.93)}
      \put(1.5,0.4){\qbezier(0,0)(3.25,1.3)(6.5,2.6)}
      \put(6.5,0.4){\qbezier(0,0)(-3.25,1.3)(-6.5,2.6)}
      \put(2.1,5){$\mathbf{7}$}
      \put(5.6,5){$\mathbf{7}$}
      \put(2.7,3.5){$\mathbf{8}$}
      \put(4.8,3.5){$\mathbf{8}$}
      \put(0.9,1.7){$\mathbf{7}$}
      \put(6.4,1.7){$\mathbf{7}$}
    \end{picture}
  \end{center}
  \caption{Smallest integral hexagon.}
  \label{min_hexagon}
\end{figure}

Besides  blowing up points to regular simplices, another technique to construct integral point sets of arbitrary dimension is to truncate simplices. By truncating regular $\di$-dimensional simplices of side length $a$ at all vertices of a regular $\di$-dimensional simplex of side length $b+2a$, we get a point set $\Ps$ with $\di^2+\di$ points. For $\di=2$, we can easily determine the set of distances of $\Ps$ to be $\{a,b,a+b,\sqrt{a^2+ab+b^2}\}$, so $\diam(\Ps)=a+b$. The smallest 
integral example is depicted in Figure \ref{min_hexagon} (here the two missing lines have edge length $7$). It is indeed the smallest integral point set with $\di=2$ and $n=6$. 

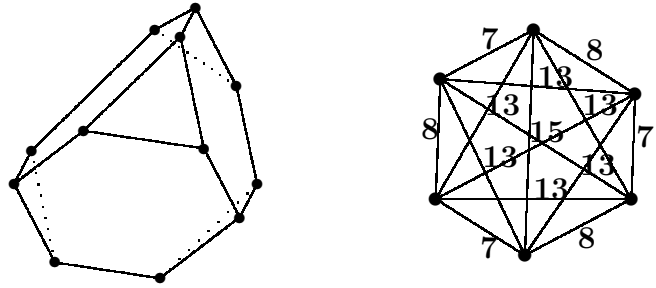
\begin{figure}[ht]
  \begin{center}
  \begin{minipage}{5.5cm}
    \setlength{\unitlength}{0.38cm}
    \begin{picture}(8.5,9.5)
      \put(0.00,3.30){\circle*{0.4}}
      \put(8.50,3.30){\circle*{0.4}}
      \put(7.78,6.72){\circle*{0.4}}
      \put(6.35,9.47){\circle*{0.4}}
      \put(7.89,2.13){\circle*{0.4}}
      \put(5.82,8.45){\circle*{0.4}}
      \put(6.64,4.54){\circle*{0.4}}
      \put(4.93,8.70){\circle*{0.4}}
      \put(5.11,0.00){\circle*{0.4}}
      \put(0.61,4.46){\circle*{0.4}}
      \put(1.43,0.55){\circle*{0.4}}
      \put(2.43,5.16){\circle*{0.4}}
      %
      \qbezier(0.00,3.30)(0.30,3.88)(0.61,4.46) 
      \qbezier(0.00,3.30)(0.72,1.92)(1.43,0.55) 
      \qbezier(0.00,3.30)(1.22,4.23)(2.43,5.16) 
      \qbezier(8.50,3.30)(8.14,5.01)(7.78,6.72) 
      \qbezier(8.50,3.30)(8.20,2.71)(7.89,2.13) 
      \qbezier[10](8.50,3.30)(6.81,1.65)(5.11,0.00) 
      \qbezier(7.78,6.72)(7.06,8.09)(6.35,9.47) 
      \qbezier[10](7.78,6.72)(6.35,7.71)(4.93,8.70) 
      \qbezier(6.35,9.47)(6.08,8.96)(5.82,8.45) 
      \qbezier(6.35,9.47)(5.64,9.09)(4.93,8.70) 
      \qbezier(7.89,2.13)(7.27,3.33)(6.64,4.54) 
      \qbezier(7.89,2.13)(6.50,1.07)(5.11,0.00) 
      \qbezier(5.82,8.45)(6.23,6.49)(6.64,4.54) 
      \qbezier(5.82,8.45)(4.12,6.80)(2.43,5.16) 
      \qbezier(6.64,4.54)(4.54,4.85)(2.43,5.16) 
      \qbezier(4.93,8.70)(2.77,6.58)(0.61,4.46) 
      \qbezier(5.11,0.00)(3.27,0.27)(1.43,0.55) 
      \qbezier[10](0.61,4.46)(1.02,2.50)(1.43,0.55) 
    \end{picture}
  \end{minipage}
  \begin{minipage}{2.7cm}
    \setlength{\unitlength}{0.40cm}
    \begin{picture}(6.7,7.5)
      \put(0.00,1.87){\circle*{0.4}}
      \put(6.50,1.87){\circle*{0.4}}
      \put(3.25,7.49){\circle*{0.4}}
      \put(0.15,5.86){\circle*{0.4}}
      \put(6.63,5.36){\circle*{0.4}}
      \put(2.96,0.00){\circle*{0.4}}
      \qbezier(0.00,1.87)(3.25,1.87)(6.50,1.87) 
      \qbezier(0.00,1.87)(1.63,4.68)(3.25,7.49) 
      \qbezier(0.00,1.87)(0.08,3.86)(0.15,5.86) 
      \qbezier(0.00,1.87)(3.32,3.61)(6.63,5.36) 
      \qbezier(0.00,1.87)(1.48,0.93)(2.96,0.00) 
      \qbezier(6.50,1.87)(4.88,4.68)(3.25,7.49) 
      \qbezier(6.50,1.87)(3.33,3.86)(0.15,5.86) 
      \qbezier(6.50,1.87)(6.57,3.61)(6.63,5.36) 
      \qbezier(6.50,1.87)(4.73,0.93)(2.96,0.00) 
      \qbezier(3.25,7.49)(1.70,6.68)(0.15,5.86) 
      \qbezier(3.25,7.49)(4.94,6.43)(6.63,5.36) 
      \qbezier(3.25,7.49)(3.11,3.75)(2.96,0.00) 
      \qbezier(0.15,5.86)(3.39,5.61)(6.63,5.36) 
      \qbezier(0.15,5.86)(1.56,2.93)(2.96,0.00) 
      \qbezier(6.63,5.36)(4.80,2.68)(2.96,0.00) 
      \put(3.25,1.87){$\mathbf{13}$} 
      \put(1.63,4.68){$\mathbf{13}$} 
      \put(-.48,3.86){$\mathbf{8}$}  
      \put(1.48,-0.09){$\mathbf{7}$}  
      \put(4.88,4.68){$\mathbf{13}$} 
      \put(6.67,3.61){$\mathbf{7}$}  
      \put(4.73,0.23){$\mathbf{8}$}  
      \put(1.50,6.88){$\mathbf{7}$}  
      \put(4.99,6.48){$\mathbf{8}$}  
      \put(3.11,3.75){$\mathbf{15}$} 
      \put(3.39,5.61){$\mathbf{13}$} 
      \put(1.56,2.93){$\mathbf{13}$} 
      \put(4.80,2.68){$\mathbf{13}$} 
    \end{picture}
  \end{minipage}
  \caption{$3$-dimensional integral point sets from a truncated tetrahedron.}
  \label{fig_3_12_min}
  \end{center}
\end{figure}

For $\di\ge 3$, the occurring distances of $\Ps$ are given by $$\mathcal{D}=\{a,b,a+b,\sqrt{a^2+ab+b^2}, \sqrt{a^2+2ab+2b^2}\},$$ so $\diam(\Ps)=\sqrt{a^2+2ab+2b^2}$. The smallest integral solution is given by $a=7$ and $b=8$ which lead to the  $\di$-dimensional integral point set with diameter $17$ consisting of $\di^2+\di$ points and proves Theorem \ref{thm_core}(a). We have depicted this integral point set for $\di=2$ and $\di=3$ in Figure \ref{fig_3_12_min}.

$\{7,8\}$ and $\{2\,021\,231,8\,109\,409\}$ are the only coprime pairs of integers with $a,b\le 10\,000\,000$, where 
all values of $\mathcal{D}$ are integers. It is not known whether infinitely many such parameter sets exist. We remark 
that a generalization of this approach to the other platonic solids does not lead to integral point sets by our methods.

\begin{figure}[ht]
  \begin{center}
    \setlength{\unitlength}{1.00cm}
    \begin{picture}(9.7,4.4)
      \put(0.00000,0.00000){\circle*{0.2}}
      \put(9.62500,0.00000){\circle*{0.2}}
      \put(5.37500,3.87298){\circle*{0.2}}
      \put(4.37500,0.00000){\circle*{0.2}}
      \put(5.25000,0.00000){\circle*{0.2}}
      \put(1.12500,0.00000){\circle*{0.2}}
      \put(2.62500,0.00000){\circle*{0.2}}
      \put(3.62500,0.00000){\circle*{0.2}}
      \put(5.50000,0.00000){\circle*{0.2}}
      \put(6.37500,0.00000){\circle*{0.2}}
      \put(7.12500,0.00000){\circle*{0.2}}
      \put(8.12500,0.00000){\circle*{0.2}}
      \put(-0.4,0){$\mathbf{A}$}
      \put(9.73,0){$\mathbf{B}$}
      \put(5.48500,3.87298){$\mathbf{P_i}$}
      \qbezier(0.00000,0.00000)(4.81250,0.00000)(9.62500,0.00000) 
      \qbezier(0.00000,0.00000)(2.68750,1.93649)(5.37500,3.87298) 
      \qbezier(9.62500,0.00000)(7.50000,1.93649)(5.37500,3.87298) 
      \qbezier(5.37500,3.87298)(4.87500,1.93649)(4.37500,0.00000) 
      \qbezier(5.37500,3.87298)(5.31250,1.93649)(5.25000,0.00000) 
      \qbezier(5.37500,3.87298)(3.25000,1.93649)(1.12500,0.00000) 
      \qbezier(5.37500,3.87298)(4.00000,1.93649)(2.62500,0.00000) 
      \qbezier(5.37500,3.87298)(4.50000,1.93649)(3.62500,0.00000) 
      \qbezier(5.37500,3.87298)(5.43750,1.93649)(5.50000,0.00000) 
      \qbezier(5.37500,3.87298)(5.87500,1.93649)(6.37500,0.00000) 
      \qbezier(5.37500,3.87298)(6.25000,1.93649)(7.12500,0.00000) 
      \qbezier(5.37500,3.87298)(6.75000,1.93649)(8.12500,0.00000) 
      \put(1.40750,1.33649){$53$} 
      \put(0.56250,0.05000){$9$}  
      \put(8.17000,1.33649){$46$} 
      \put(8.84500,0.05000){$12$} 
      \put(4.37500,1.33649){$32$} 
      \put(4.88250,1.33649){$31$} 
      \put(2.22000,1.33649){$46$} 
      \put(3.15000,1.33649){$38$} 
      \put(3.80000,1.33649){$34$} 
      \put(5.47750,1.33649){$31$} 
      \put(6.03000,1.33649){$32$} 
      \put(6.55000,1.33649){$34$} 
      \put(7.20000,1.33649){$38$} 
      \put(4.81250,0.05000){$7$}  
      \put(4.00000,0.05000){$6$}  
      \put(5.30500,0.15000){$2$}  
      \put(1.87500,0.05000){$12$} 
      \put(3.12500,0.05000){$8$}  
      \put(5.93750,0.05000){$7$}  
      \put(6.75000,0.05000){$6$}  
      \put(7.62500,0.05000){$8$}  
    \end{picture}
    \caption{$2$-dimensional integral point set with  $n=12$ and diameter $77$.}
    \label{fig_3_15_min}
  \end{center}
\end{figure}
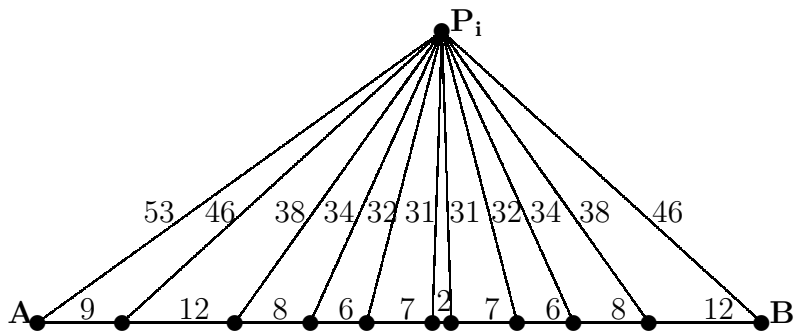

There is another important construction of integral point sets for. In Figure \ref{fig_3_15_min}, we have depicted a plane integral point set consisting of $12$ points of which $11$ are col-

\begin{figure}[ht]
  \begin{center}
    \setlength{\unitlength}{0.675cm}
    \begin{picture}(2.5,1)
      \put(0,0.0){\line(1,0){2}}
      \put(0,0.0){\circle*{0.3}}
      \put(2,0.0){\circle*{0.3}}
      \put(0,0.0){\qbezier(0,0)(0.125,.475)(0.25,.97)}
      \put(2,0.0){\qbezier(0,0)(-0.125,.475)(-0.25,.97)}
      \put(0.25,0.94){\line(1,0){1.5}}
      \put(0.25,0.93){\circle*{0.3}}
      \put(1.75,0.93){\circle*{0.3}}
      \put(0,0.0){\qbezier(0,0)(0.875,.475)(1.75,.97)}
      \put(2,0.0){\qbezier(0,0)(-0.875,.475)(-1.75,.97)}
      \put(0.9,-0.4){$4$}
      \put(-0.16,0.3){$2$}
      \put(1.94,0.3){$2$}
      \put(0.9,1.05){$3$}
      \put(0.34,0.35){$4$}
      \put(1.34,0.35){$4$}
    \end{picture}
    \caption{Smallest plane integral point set with $n=4$ and no three points on a line.}
    \label{fig_min_n3ol}
  \end{center}
\end{figure}
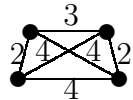

linear. If we rotate the point $P_i$ around the line $\overline{AB}$, it forms a circle 
with radius $8\sqrt{15}$. If we scale the plane integral point set of Figure \ref{fig_min_n3ol} by a factor of 
$15$, we can arrange $4$ points on this circle, so that we get a $3$-dimensional integral point set consisting of 
$15$ points with diameter $d(3,15)=77$. 

Generally, for $\di\ge 2$ we may consider an $(\di-1)$-dimensional sphere $S$ that intersects $P_i$, has its center on the line $\overline{AB}$, and spans a hyperplane that is orthogonal to $\overline{AB}$. (If $\di=2$, then $S$ consists of the point $P_i$ and its reflection in $\overline{AB}$.) If it is possible to place $n'$ points on this sphere with pairwise integral distances, then these points together with the points of the line $\overline{AB}$ form an $\di$-dimensional integral point set consisting of $n+n'-1$ points. This gives the proof of Theorem \ref{thm_line_circle}. It is tight for $\di=3$, $13\le n\le 21$. Nevertheless, we conjecture that Theorem \ref{theo_blow_up_n_minus_1} yields better bounds for $\di=3$ and $n\ge 22$.

\section*{Acknowledgement}
The authors thank the anonymous referee for his helpful comments and suggestions on the paper. Moreover, we would like to thank Nikolai Avdeev 
for pointing out the missing condition $v<2r$ in Theorem~\ref{theo_blow_up_2} to us.
 

\bibliographystyle{abbrv}
\bibliography{upper_bounds_for_integral_point_sets}
\bibdata{upper_bounds_for_integral_point_sets}

\end{document}